\def \version {2026 -- 09 -- 21}

\documentclass[12pt]{article}
\usepackage{amsthm,amssymb,amsmath}
\usepackage{lineno}
\usepackage{comment}

\def \bpf {\begin{proof}}
\def \epf {\end{proof}}
\newtheorem{theorem}{Theorem}
\def \btm {\begin{theorem}}
\def \etm {\end{theorem}}
\newtheorem{proposition}[theorem]{Proposition}
\def \bpn {\begin{proposition}}
\def \epn {\end{proposition}}
\newtheorem{corollary}[theorem]{Corollary}
\def \bcr {\begin{corollary}}
\def \ecr {\end{corollary}}
\newtheorem{lemma}[theorem]{Lemma}
\def \blm {\begin{lemma}}
\def \elm {\end{lemma}}
\newtheorem{problem}[theorem]{Problem}
\def \bpm {\begin{problem}}
\def \epm {\end{problem}}
\newtheorem{definition}[theorem]{Definition}
\def \bdf {\begin{definition}}
\def \edf {\end{definition}}
\newtheorem{remark}[theorem]{Remark}
\def \brm {\begin{remark}}
\def \erm {\end{remark}}
\newtheorem{example}[theorem]{Example}
\def \bex {\begin{example}\rm }
\def \eex {\end{example}}
\newcommand{\nev}[1]{{\bf\itshape (#1)} \ }
\newcommand{\theoa}[4]{\nin \textbf{Theorem #1} \textrm{(#2)} \nev{#3} \textit{#4}}
\newcommand{\theos}[3]{\nin \textbf{Theorem #1} \textrm{(#2)} \textit{#3}}

\def \bsk {\bigskip}
\def \msk {\medskip}
\def \ssk {\smallskip}
\def \igs {\bsk}
\def \nin {\noindent}
\def \smin {\diagdown}	

\def \cF {\mathcal{F}}

\def \st {\mathrm{st}}
\def \ex {\mathrm{ex}}
\newcommand{\floor}[1]{\lfloor #1 \rfloor}
\newcommand{\ceil}[1]{\lceil #1 \rceil}

\title{\vspace*{-2.5cm} ~~~ \\
The strong (non-induced) Tur\'an numbers}
\author{Yair Caro$^1$, Zsolt Tuza$^{2,3}$ \\
~~~ \\
\normalsize $^1$ Department of Mathematics, \vspace*{1ex} University of Haifa-Oranim, Israel \\ 
\normalsize $^2$ HUN-REN \vspace*{1ex} Alfr\'ed R\'enyi Institute of Mathematics, Budapest, Hungary \\
\normalsize $^3$ Department of Computer Science and Systems Technology \\ \normalsize \vspace*{1ex} University of Pannonia, Hungary}

\date{\small Latest update on \version }

\begin{document}

\maketitle

\begin{abstract}
In this paper we introduce and explore the following new graph invariant:
For a graph $G$ on $k$ vertices, $G \neq  K_k$, let $\st(n,G)$
 denote the maximum number of edges in a graph of order $n$ which
 does not contain any subgraph on $k$ vertices strictly containing $G$.

A basic relation to classical Tur\'an numbers is developed via the following:
For $G$ on $k$ vertices, let $D(G) = \{ H : |H| = |G|, \, H = G + e \}$.   
Using this notion we prove that
 $\ex(n,G) \leq \st(n,G) = \ex(n, D(G) ) \leq \min \{ \ex(n,H) : H \in D(G) \}$ holds for all $n \geq |G|$.   

The family $D(G)$ happened to be smoothly amenable to the use of classical extremal results, and in many cases allows us to get asymptotically sharp estimates as well as exact values of $\st(n,G)$.

From the many results proved here we state the following as an illustration.

(1)\quad If $\chi(G) \geq 3$ and $\chi(D(G)) = \chi(G)$, then $\st(n,G) =   (1+o(1))\,\ex(n,K_{\chi(G)})$.  

(2)\quad If $\chi(D(G)) = \chi(G) +1$, then $G$ is a complete $\chi(G)$-partite graph and $\st(n,G) = \ex(n,K_{\chi(G) +1})$ for $n$ sufficiently large.  

(3)\quad For $k$ odd, $k\geq 5$, $\st(n,C_k) = \ex(n,C_k)  = \ex(n,K_3)$ for $n$ sufficiently large.  

(4)\quad If $T$ is a tree of order $q$ with diameter $k \geq 2$ and $q \geq k+1 \geq 3$,  then $\ex(n, \{C_3,\dots,C_{k+1}) \leq \st(n,T) \leq \ex(n, \{C_3,\dots,C_{k+1})  + (q-1)n$.

Many results concerning even cycles, theta graphs,  dense bipartite graphs and graphs of the form $G = G^* \cup tK_1$ are obtained, moreover the value of $\st(n,G)$ is computed for all graphs on at most 4 vertices.
\end{abstract}

\section{Introduction}
\label{s:intro}

We consider the following extremal Tur\'an-type problem.
 
For a graph $G \neq  K_k$ on $k$ vertices, and for every integer $n \geq  k$, let $\st(n,G)$ be the largest integer $m$ for which there exists a graph $H$ with $n$ vertices and $m$ edges such that no 
 $F \subset H$ on $k$ vertices contains $G$ as a \emph{proper} subgraph.
So, $\st(n,G)$ measures whether and how much we can pass $\ex(n,G)  = \max \{ m : \exists H, \, |H|=n, \, e(H) = m, \, G \not\subset H \}$,
without the appearance of a subgraph $F$ on $|G|$  vertices strictly containing $G$ but allowing that $G$ itself is present.

Classical Tur\'an theory asks how many edges are possible while avoiding $G$, whereas the
 strong Tur\'an number asks how many edges are possible while allowing $G$, but forbidding every graph obtained from $G$ by adding an edge.
Surprisingly, this apparently small relaxation can radically change the extremal behaviour, and the resulting parameter is governed by a family denoted $D(G)$, which we will introduce formally soon.

A major distinction between $\ex(n,G)$ and $\st(n,G)$ is monotonicity.
While $\ex(n,G)$ is monotone, namely $G\subset H$ implies $\ex(n,G) \leq  \ex(n,H)$, it is
far from being the case for $\st(n,G)$.
A bold example is given by the path $P_3$, for which $\st(n,P_3) = \ex(n,K_3) = \floor {n^2\!/4}$ while $\st(n, P_3 \cup K_1)  = \floor {2n/3}$.
Another distinction from classical Tur\'an problems is raised by the role of isolated vertices.
While on the classical Tur\'an function isolated vertices have no impact, namely $\ex(n,G) = \ex(n, G \cup tK_1)$ for all $n \geq   |G| +t$, isolated vertices may have substantial impact
on $\st(n,G)$ as already demonstrated in $\st(n,P_3)$ vs.\ $\st(n,P_3 \cup K_1)$.

A basic relation connecting $\st(n,G)$ to classical Tur\'an  problems is obtained via the following definition.
For $G \neq  K_k$, of order $|G| = k$, let $D(G) = \{H : |H| = |G|, \, H = G + e\}$. 
A simple but very useful equivalence, proved as
 Lemma \ref{l:equiv} in Section \ref{s:general},
states that $\st(n,G) = \ex(n,D(G))$.

The structure of $D(G)$ makes it amenable to the use of many classical
 results in external graph theory, among them the
  Erd\H os--Stone--Simonovits theorem \cite{ESto,ESim},
Simonovits' critical edge theorem \cite{S68}
 and Dirac's theorem \cite{D63}.
The current work is also inspired in part by a recent paper
  ``Induced / Incomparable versus Ramsey'' \cite{CTZ}
where similar problems are treated in the context of Ramsey theory.

\subsection{Organization of the paper}

In Section \ref{s:general} we give general results obtained
 using classical extremal theorems mentioned above, and in particular
  among others, we prove items (1) and (2) mentioned in the Abstract.

Section \ref{s:isol} is devoted to the role and impact of isolated vertices.
We completely determine $\st(n , G \cup  tK_1)$ in terms of the Tur\'an function for $t \geq  2$, and determine $\st(n,G \cup K_1)$ up to an additive linear function of $n$.

Section \ref{s:trees} deals with forests.
We explicitly determine $\st(n,G)$ for $G =  K_{1,k}$, $tK_2$, $P_4$, and prove the following asymptotic sharp theorem:  
Let $T$ be a tree with diameter $k \geq  2$ and order $q \geq  k +1 \geq  3$.
Then $\ex(n, \{C_3, \dots ,C_{k+1}\}) \leq  \st(n, T) \leq  \ex(n, \{C_3, \dots ,C_{k+1}\}) + (q - 1)n$.

Section \ref{s:theta} is devoted to the study of $\st(n,C_{2k})$ for $k \geq  3$  (the case of $C_4$ is solved exactly in Section \ref{s:general}).
Using a recent breakthrough  concerning theta graphs \cite{LY23}
and close consideration of $D(C_{2k})$, we determine the correct order of magnitude of $\st(n,C_{2k})$ for $k =  3,4,5, 6, 9,10$
and otherwise lower and upper bounds not too far apart.

Section \ref{s:bip} concerns applications of an old theorem of Erd\H os
 \cite{Er65} to approximate $\st(n,G)$ when $G$ is sandwiched between $K_{p, p}$ and $K_{p+1,p+1}$ from which two edges are deleted.

In section \ref{s:open} we present one example of many possible directions for future research.

\subsection{Notation}

We shall use standard graph theory notation \cite{west}.
In particular, the number of vertices (order of $G$) and number of edges are denoted as $|G|$ and $e(G)$, respectively.
Furthermore, $\deg(v)$ denotes the degree of vertex $v$, and $\delta = \delta(G)$, $\Delta = \Delta(G)$ and $d = d(G) = 2e(G)/|G|$ are respectively the minimum degree, the maximum degree, and the average degree of $G$.
Other notation will be explained when first appears in the paper.

\newpage

\section{General bounds from classical extremal graph theory}
\label{s:general}

In this section we derive bounds on $\st(n,G)$ in terms of classical extremal graph theory.

Recall the important notation
 $D(G) = \{H : |H| = |G|, \, H = G+e\}$ for $G \neq K_k$; namely, $D(G)$ is the collection of graphs having the same order as $G$, obtained from $G$ by inserting any one new edge.
Define $\chi(D(G)) = \min\{\chi(H) : H \in D(G)\}$.
 
Our first result connects $\st(n,G)$ with $\ex(n,D(G))$.

\blm
\nev{Equivalence Lemma}
\label{l:equiv}
For every graph $G$ we have $\ex(n,G) \leq  \st(n,G) = \ex(n, D(G) ) \leq   \min \{ \ex(n,H) : H \in  D(G) \}$. 
\elm

\bpf
Lower bound: 
$\ex(n,D(G)$  means that there exists a graph $H$  with $n$ vertices and $\ex(n,D(G))$ edges with no copy of any member of $D(G)$, hence no non-induced subgraph on $k$ vertices contains $G$ in $H$. 

Upper bound: 
$\ex(n,D(G))+1$ means that every graph on $n$ vertices and $\ex(n,D(G))+1$ edges contains a (not necessarily induced) copy of a member of $D(G)$, hence contains a non-induced copy of $G$. 
The rightmost inequality is trivial. 
\epf

So, computing $\st(n,G)$ is equivalent to computing $\ex(n,D(G))$, however when bipartite graphs are involved in $D(G)$ computing $D(G)$ even asymptotically can be a hard task, in many cases unknown yet.  

The following classical theorems will be used frequently in the sequel. 

\igs

\theoa{A.}{Erd\H os, Stone \cite{ESto}; Erd\H os, Simonovits \cite{ESim}}{Erd\H os--Stone--Simo\-no\-vits Theorem, ESS}{Let $\cF = \{ G_1,\dots,G_k \}$ be a family of graphs with $\chi(\cF) := \min \{\chi(G_ j) \mid j  = 1,\dots,k \} = r  \geq 3$.  Then $\ex(n,\cF) = (1+o(1))\,\ex(n,K_r ) = \frac{r-2}{2r-2}\,n + o(n^2)$.}

\igs

An edge $e$ in a graph $G$ is called color-critical if $\chi(G-e)<\chi(G)$.

\igs

\theoa{B.}{Simonovits \cite{S68}}{Simonovits' Color-Critical Edge Theorem, SCCE}{If $H$ contains a color-critical edge and $\chi(H) = k+1$, then there exists an $n_0(H)$ such that $\ex(n, H) = e(T_{n,k})$ and the $k$-partite Tur\'an graph $T_{n,k}$ is the only extremal graph, provided $n \geq n_0(H)$.}

\igs

We also recall the following result which, although less general, is tight for all $n>r$.

\igs

\theos{C.}{Dirac \cite{D63}}{For every $r\geq 3$ and $n \geq  r +1$, every graph with $n$ vertices and $\ex(n,K_r) +1$ edges contains $K_{r+1} - e$.}

\igs

We are now ready to present the main theorem of this section.

\btm
~

\begin{itemize}
 \item[$(i)$] If $\chi(G)  \geq  3$ and $\chi(D(G)) =  \chi(G)$, then $\st(n,G) = (1+o(1))\ex(n,K_{\chi(G)})$.
 \item[$(ii)$] If  $\chi(D(G))   = \chi(G) +1$,  then  $G$ is a complete $\chi(G)$-partite graph and  $\st(n,G) = \ex(n,K_{\chi(G) +1})$ for every $n$  sufficiently large.
 \item[$(iii)$] If $\chi(G)= 2$ and $\chi(D(G)) = 2$, then $\st(n,G) \leq  \ex(n,H)$ where $H$ is complete bipartite containing $G$ as a spanning non-induced subgraph, and $\st(n,G)$ is sub-quadratic.
\end{itemize}
\etm

\bpf

$(i)$\quad 
By Theorem A  and Lemma \ref{l:equiv} we obtain $(1+o(1))\ex(n,K_{\chi (G)}) = \ex(n,G ) \leq  \st(n,G)  =  \ex(n,D(G))  =  (1+o(1))\ex(n,K_{\chi (G)})$  as $\chi (G) = \chi (D(G))$.

$(ii)$\quad 
If $\chi(G) = 1$, then $G$  is the empty graph and the result is immediate. 
Otherwise, $\chi(G)\geq 2$.
If $G$ is not a complete $\chi(G)$-partite graph with $\chi(G) \geq  2$, then embed it in a complete $\chi(G)$-partite $H$ with $|H| = |G|$  (using the $\chi(G)$ color classes of $G$ from a proper coloring as parts).
There is a missing edge in $H$ that can be added to get a graph $G^*$ with $|G^*| = |G|$, $G^*$ strictly containing $G$ and with $\chi(G^*) = \chi(G)$, contradicting the assumption $\chi(D(G))  = \chi(G) +1$.
On the other hand, if $G$ is a complete $\chi(G)$-partite graph, then adding any edge is possible only within the color classes, hence forcing an increase of the chromatic number by~1. 
Consequently, all these graphs formed  by adding one edge contain a color-critical edge, namely the deletion of the added edge reduces the chromatic number by 1 back to $\chi(G)$. 
The rest follows from Theorem B applied for a family of critical graphs instead of just one graph. 

$(iii)$\quad 
According to the assumptions,  $G$ is bipartite but not complete bipartite; hence it can be extended to a complete bipartite $H$ by adding the edges missing between the two color classes of $G$.  
Now, by definition, $\st(n,G)  \leq  \ex(n,H)$ holds.
The sub-quadratic upper bound on $\st(n,G)$ follows by the theorem of K\H ov\'ari, S\'os and Tur\'an \cite{KST}.
\epf

We next determine $\st(n,C_k)$ for odd $k \geq 5$ and $n$ large enough.

\bpn
For $k$  odd, $k \geq  5$ and for $n$ large enough,  $\st(n,C_k) = \ex(n,C_k)  = \ex(n,K_3)$.
\epn

\bpf  
Observe that for an odd $k$-cycle ($k \geq  5$) adding any chord $e$ makes $C_k$ contain a copy of a shorter odd cycle $C^*$  and an even cycle sharing the common edge $e$.
Deleting any edge other than $e$ from the odd cycle $C^*$ leaves a bipartite graph. 
Hence,  every edge of  $C^*$ except the chord is critical and by SCCE  for all odd $k \geq  5$ and $n$ large $\st(n,C_k) = \ex(n,C_k)  = \ex(n,K_3)$ holds.
\epf

Dirac's  theorem enables us to determine exact values of $\st(n,G)$ for further graphs $G$, as follows.

\bpn
\label{p:st=ex}
~

\begin{itemize}
 \item[$(i)$] Suppose $G  =  K_{r+1} - 2K_2$   or  $G = K_{r +1} - P_3$.
Then $\st(n,G)  = \ex(n,K_r)$ for all $n \geq r+1$.
 \item[$(ii)$] Suppose $r \geq 2$, and let $G$ be a graph on $r+1$ vertices such that $\chi(G) = r$ and $e(G) \leq \binom{r+1}{2} - 2$.  Then $\st(n,G)  = \ex(n,K_r)$ holds for all $n \geq r +1$.
\end{itemize}
\epn

\bpf
Observe that in both cases of $(i)$, $K_{r+1} - e$ is the unique member of D(G). Applying the Equivalence Lemma we obtain $\st(n,G)  =  \ex(n, K_{r+1} - e)   = \ex(n,K_r)$ by Theorem C. 

For $(ii)$ we observe that all graphs in $D(G)$ are on $r+1$ vertices and having chromatic number $r$, as they contain the $r$-chromatic $G$ but are strictly  contained in $K_{r+1}$.
Hence the combination of Tur\'an's theorem, Theorem~C, and Lemma \ref{l:equiv} yields $\ex(n,K_r) \leq  \ex(n,G) \leq  \st(n,G)  =  \ex(n,D(G)) =  \ex(n,K_r )$.
\epf

\newpage

\section{The role/impact of isolated vertices}
\label{s:isol}

In this section we analyze the behavior of $\st(n,G)$ on graphs which
 contain isolated vertices.
As a triviality, one can see that $\st(n,E_k) = 0$ holds for all $n\geq k > 0$.
We next give some simple observations, which will turn out to be useful concerning $\st(n,G)$.

\blm
\label{l:diff}
For every $m$ in the range $2\leq m\leq\binom{k}{2}-2$ there exist
 at least two non-isomorphic graphs with $k$ vertices and $m$ edges.
\elm

\bpf
The considered range implies $k\geq 4$.
Hence for $m=2$ the graphs $P_3$ and $2K_2$ supplemented with
 $k-3$ or $k-4$ isolated vertices do the job, respectively.
If $3\leq m\leq \frac{1}{2}\binom{k}{2}$, $m$-edge triangle-free graphs
 (e.g., bipartite subgraphs of $K_{\floor{n/2},\ceil{n/2}}$) and also $m$-edge graphs containing $K_3$ exist;
 clearly the two cannot be isomorphic.
Finally, for any larger $m$ we can take the complementary graphs of the above.
\epf

As an immediate consequence, for every $G$ with $|G|=k$ and $2\leq e(G)\leq \binom{k}{2}-2$
  we have $\ex(n,G)\geq e(G)$.
If $G$ is assumed to be connected, then the following stronger assertion can also be proved.

\blm
\label{l:ex>e2}
If $G$ is a connected graph of order $|G| = k \geq  4$, then
 $\ex(n,G) \geq   \ex(n - k,G) + e(G)$ for all $n \geq  k$,
unless $n=k$ and $G=K_k$ or $G=K_k-e$.
\elm

\bpf
Assume first $G\notin \{ K_k , K_k-e \}$.
Since $G$ is connected, it has at least $k-1\geq 3$ edges, hence the
 conditions of Lemma \ref{l:diff} are satisfied and there exists
 a graph $H\ncong G$ with $|H|=|G|$ and $e(H)=e(G)$.
Take a graph $F$ realizing $\ex(n - k,G)$ (even when $n - k = 0, 1, 2,\dots$ is small) and add $H$ on $k$ further vertices.
As $G$ is connected, no copy of $G$ can emerge, thus $\ex(n,G) \geq  \ex(n- k, G) + e(G)$.

If $G\in \{ K_k , K_k-e \}$, let $H$ be any 2-connected graph
 obtained from $G$ by the removal of one edge.
If $G=K_4-e$, we have $H=C_4$, and in all the other cases
 any one edge can be deleted from $G$ to get $H$.
As previously, also now we take $F\cup H$, additionally joining
 $F$ and $H$ with a further edge.
No copy of $G$ can emerge because $G$ is 2-connected, hence
 cannot contain the last inserted cut-edge.
\epf

\newpage

\subsection{More than one isolated vertex: \textbf{\itshape G}\/\,$=$\,\textbf{\itshape G}\/$\mathbf{^*}$\,$\cup$\,\textbf{\itshape tK}$\mathbf{_1}$, \textbf{\itshape t}\,$\geq$\textbf{2}}

\btm
Let $G = G^*\cup tK_1$ where $t \geq 2$ and $G^*$ is a graph with at least two edges.
Then $\st(n,G)  = \ex(n,G^*)$ holds for all $n\geq |G|$.
\etm

\bpf
Since every $G^*$-free graph is $G$-free, we always have $\st(n,G) \geq \ex(n,G^*)$.
To derive the converse inequality, consider any graph $H$ of order
 $n\geq |G|$ with more than $\ex(n,G^*)$ edges.
We prove that $H$ has a subgraph $F$ of order $|G|$ which
 strictly contains $G$, hence $\st(n,G)<\ex(n,G^*)+1\leq e(H)$ will follow, implying the theorem.

For the proof we first observe that the conditions $e(G^*)\geq 2$ and
 $n\geq |G| \geq |G^*|+2$ imply $\ex(n,G^*) \geq e(G^*)$.
Indeed, if $G^*$ is not a matching, then we can pick a non-isolated edge
 say $e$ from $E(G^*)$ and construct the graph $(G^*-e)\cup K_2\cup (n-2-|G^*|)K_1$;
  it has $e(G^*)$ edges and is $G$-free (having more isolated edges than $G^*$).
And if $G^*=sK_2$ is a matching of size $s\geq 2$, then 
 $K_{2s-1}\cup (n-2s+1)K_1$ is $G$-free and has $\binom{2s-1}{2}>s$ edges.

Since $e(H)>\ex(n,G^*)\geq e(G^*)$, we can find a subgraph $G'$ isomorphic to
 $G^*$ in $H$, moreover $H$ contains a further edge say $e'$ in
 $E(H)\smin E(G')$.
Let $e' = xy$  and $k = 2  -  |\{x ,y \} \cap V(G')|$.
Then, supplementing $G'\cup e'$ with $|G|-|G'|-k$ vertices from
 $V(H)\smin V(G'\cup e')$ we obtain the required subgraph $F$ of $H$,
 which strictly contains $G$.
\epf

\brm
For the leftover case $G^*=K_2$, i.e.\
 $G = K_2 \cup tK_1$, $t \geq 2$, we have $|G|= t +2 \geq 4$, hence every graph with $n \geq |G|$ vertices and at least two edges must  contain $P_3 \cup (t-1)K_1$ or $2K_2 \cup (t-2)K_1$,
both of them strictly containing $G$.
Therefore, $\st(n, K_2 \cup tK_1)  = 1  = \ex(n, G^*) +1$ holds for all $t \geq 2$ and $n \geq t +2$.
\erm

\subsection{One isolated vertex: \textbf{\itshape G}\/\,$=$\,\textbf{\itshape G}\/$\mathbf{^*}$\,$\cup$\,\textbf{\itshape K}$\mathbf{_1}$}

\subsubsection*{A --- The connected case of \textbf{\itshape G}\/$\mathbf{^*}$}

\btm
Suppose  $G = G^*\cup K_1$, where $G^*$ is connected, and let $n\geq |G|$.
 \begin{itemize}
  \item[$(i)$] If $|G^*| \geq  4$,  then for all $n \geq  n_0(G^*)$ we have $\st(n,G) = \ex(n,G^*)$.
  \item[$(ii)$] If $G^* = K_2$, then $\st(n,K_2 \cup K_1) = \ex(n, P_3) = \floor {n/2}$.
  \item[$(iii)$] If $G^* = P_3$,  then $\st(n,P_3 \cup K_1)) = \floor{2n/3}$.
  \item[$(iv)$] If $G^* = K_3$, then $\st(n, K_3\cup K_1) = \ex(n,K_4-P_3) = \ex(n,K_3) = \floor{n^2\!/4}$.
 \end{itemize}
\etm

\bpf
Clearly $\st(n,G)  \geq  \max \{ \ex(n,G^*)  , |G| \}$.
This fact directly proves the lower bounds in $(i)$ and $(iv)$.
For $(ii)$ and $(iii)$ one can take $\floor {n/|G^*|}$ copies of $K_2$
 or $P_3$, respectively, with a further isolated vertex or isolated edge
 if $n$ is not a multiple of $|G^*|$, i.e.\ if $2\nmid n$ or $3\nmid n$.

Turning to the corresponding upper bounds we consider the four cases separately,
 assuming $n \geq   |G |$ and considering any graph $H$ whose $e(H)-1$
 equals the claimed value of $\st(n,G)$.
We will prove that such an $H$ strictly contains $G$.
 
\msk

\nin
$(i)$\quad
Let $|G^*| \geq  4$ and $e(H) \geq  \ex(n, G^*) +1$.

\ssk

\nin
\underline{Case A:} \ $G^* \notin  \{  K_k  , K_k -e  \}$.

Then $H$  contains a copy of $G^*$ which must be induced and isolated,  namely no edge outside $E(G^*)$ is incident with a vertex of $G^*$.  Otherwise, we obviously have a non-induced copy of $G$. 
Removing $V(G^*)$ we get a graph $H^*$ on $n - k$ vertices and $e(H^*) = \ex(n, G^*) + 1 - e(G^*) > \ex(n-k,G^*)$, by Lemma \ref{l:ex>e2}.
So, $H^*$ contains another copy of isolated $G^*$.  We can continue deleting isolated copies of $G^*$  until we are left with at most  $k -1$  vertices (or else along the way we get a non-induced copy of $G$).  
Since if we have at least $t \geq   k$ vertices and at least $\ex(t ,G^*)+1$ edges,  then again, an isolated copy of $G^*$ exists and we may delete it.  Once we are left with at most $k -1$ vertices then  for  $j < k$ we can have at most $\binom{j}{2}$ edges.
So, the total number of edges in $H$  where $n =  rk +j$, $j = 0 ,\dots,k-1$, is at most  $\binom{j}{2} +  e(G^*)n/k \leq  \binom{k-1}{2} + e(G^*)n/k$.   

If $G^*$ is acyclic, then this number is at most $\binom{k-1}{2} +  n(k-1)/k$  but for all $G^*$ on at least 4 vertices---especially for $P_4$ and $K_{1,3}$---we have   $\ex(n,G^*) \geq  n -1$, a contradiction for $n > k(k-1)(k-2)/2$.

If $G^*$ is not acyclic, then it contains a cycle of length at most $k$, hence $\ex(n,G^*)$ is not linear in $n$. 
Then crudely $\binom{j}{2} + e(G^*)n/k \leq  \binom{k-1}{2} +k(k-1)n/2k = \binom{k-1}{2}  +  n(k-1)/2$, which is linear in $n$,  so it is a contradiction for $n$ large enough. 

\ssk

\nin
\underline{Case B:} \ $G^* \in  \{  K_k  , K_k -e  \}$.

Here essentially the same proof applies, except that Lemma \ref{l:ex>e2}
 makes a distinction for $n=k$, therefore
 we have to stop when at most $k +1$ vertices are reached instead of $k$. 
However here the comparison of $\binom{k}{2} + \binom{k}{2}n/k$  is against $\ex(n,G^*)$ which is quadratic, due to Theorem C.
Thus, the contradiction follows that $e(H)$ is much smaller than $\ex(n,G^*)$ if $n$ is large.
This completes the proof of $(i)$.

\msk

\nin
$(ii)$\quad
If $G^* = K_2$, $G = K_2 \cup K_1$, then $D( G) = \{P_3 \}$, and of course every $H$ with $e(H) = \ex(n, P_3)+1 = \floor{n/2}+1$ edges contains $P_3$.
Thus, $\st(n,G)\leq \ex(n,P_3)$ and as we have seen, equality must hold due to the graphs $\floor{n/2}K_2$ and $\floor{n/2}K_2\cup K_1$ depending on the parity of $n$.

\msk

\nin
$(iii)$\quad
If $G^* = P_3$, $G = P_3 \cup K_1$, then $D(G) = \{P_4, K_{1,3}, C_3\}$.
Consider any $D(G)$-free graph $H$ on $n\geq 4$ vertices.
Then each component of $H$ is either $K_1$ or $K_2$  or $P_3$. 
The largest edge/vertex ratio occurs for $P_3$, so the maximum of $e(H)$ is attained where $\floor{n/3}$ copies of $P_3$ are taken, together with $K_{n- 3\floor{n/3}}$ if $3\nmid n$. If $n\equiv 1$ (mod 3), then an alternative is to replace one $P_3$ and the remaining $K_1$ with $2K_2$.

\msk

\nin
$(iv)$\quad
If $G^* = K_3$, $G = K_3 \cup K_1$, then $D(G)$ is the paw graph $K_4-P_3$, and clearly $\st(n,G) \geq \ex(n,K_3)$.
In any $(K_4-P_3)$-free graph $H$, every triangle is isolated.
Assume $H$ has $t$ triangles; then $e(H) \leq 3t + \ex(n-3t, K_3) = 3t + \floor{(n-3t)^2\!/4} \leq \floor{n^2\!/4}$, with equality if and only if $t=0$ and $H$ is the bipartite Tur\'an graph.
\epf

\subsubsection*{B --- Unconditional on the connectivity of \textbf{\itshape G}\/$\mathbf{^*}$}

In the following we consider $G = G^* \cup  K_1$ with
 $\delta (G^*) > 0$, where no connectivity is assumed for $G^*$.

\btm
Suppose $G =  G^* \cup K_1$  where $G^*$ has $m$ edges, $k$ vertices and $\delta (G^*) > 0$.
Then $\ex(n,G^*) \leq  \st(n,G) \leq  \ex(n,G^*) + \min \{ \delta (G^*)n , nm/k \}$.   
\etm

\bpf
The lower bound is immediate. 
For the upper bound suppose $H$ has $n$ vertices and $\ex(n,G^*) + \delta (G^*)n  +1$ edges, without a non-induced copy of $G$. 
Repeatedly remove from $H$ vertices of degree at most $\delta (G^*)$ until no such vertex remains; and call the remaining graph $F$.

Since we removed at most $\delta (G^*)n$ edges, fewer than $e(H)$, the
 graph $F$ is nonempty.
Also, if $F$ contained no copy of $G ^*$,  then by monotonicity
 of the Tur\'an function we would have
 $e(F) \leq  \ex(|F|,G^*) \leq  \ex(n, G^*)$, so that
     $e(H) \leq  e(F) + \delta (G^*)n  \leq  \ex(n,G^*) + \delta (G^*)n$ would follow,
     contradicting $e(H) = \ex(n,G^*) + \delta (G^*)n +1$.

Hence, we can choose a copy of $G^*$ in $F$ and also select a
 vertex $v$ having degree $\delta (G^*)$ in this copy.
As all degrees in $F$ are at least $\delta (G^*) +1$, the selected $v$
 is incident with a further edge in $F$, going either to another vertex of $G^*$ or out to provide $G^*$ with an attached leaf.
In either case we have a non-induced copy of $G =  G^* \cup K_1$.

Next, consider $H$ with $n$ vertices and $\ex(n,G^* ) + nm/k +1$ edges. 
Then $H$ contains a copy of $G^*$, which must be induced and isolated (no edge out), otherwise a non-induced copy of $G$ would occur.
Repeatedly remove copies of $G^*$  until no more copies of $G^*$ are present; call the remaining graph $F$.

Similarly as above, $F$ cannot be empty, since we deleted $nm/k$ edges, fewer than $e(H)$.
If $F$ contains no copy of $G^*$, then $e(F) \leq  \ex(|F|,G^*) \leq  \ex(n,G^*)$ as the  Tur\'an function is monotone increasing. 
But having deleted $t \leq  n/k$ copies of $G^*$, we obtain $e(F) + tm \leq  \ex(n,G^* ) + nm/k < e(H)$, a contradiction.
Thus, $\st(n,G) \leq  \ex(n,G^*) +nm/k$.    
\epf

\section{Forests, trees and diameter}
\label{s:trees}

This section is devoted to the study of $\st(n,G)$ on various kinds of forests.
It turns out that the diameter of a tree $T$ plays an important role in estimating $\st(n,T)$.
Before proving that, we consider stars, matchings, and the path $P_4$.

\btm
For the $k$-star $K_{1,k}$ we have
 $$
   \st(n,K_{1,k}) =
     \begin{cases}
  \begin{tabular}{ccl}
    $\floor{n^2\!/4} = \ex(n,K_3)$ && \textrm{for all } $n \geq  2k-1$, \\
    $\floor{(k-1)n/2}$ && \textrm{if } $k+1\leq n \leq  2k-2$.
  \end{tabular}
     \end{cases}
 $$
\etm

\bpf
Consider any graph $H=(V,E)$ with $n$ vertices and $m$ edges.
Then summing up the degrees along the edges we obtain
 $$
   \sum_{vw\in E} (d(v)+d(w)) = \sum_{v\in V} (d(v))^2
     \geq
   n \left( \frac{\sum_{v\in V} d(v)}{n} \right) ^{\hspace*{-0.2em}2} = \frac{4m^2}{n} \,.
 $$
Thus, there exists an edge $vw$ with $d(v)+d(w) \geq \ceil{4m/n}$.

If $n \geq 2k-1$ and $m>\floor{n^2\!/4}$, this means a $vw$
 with $d(v)+d(w) \geq \ceil{4m/n} > n \geq 2k-1$.
Then $vw$ is contained in a triangle, say $vwz$.
Assuming $d(v)\geq d(w)$ we also have $d(v)\geq k$, hence $vwz$ can be
 extended to a $K_{1,k}$ with a further edge in the neighborhood of $v$.
This implies $\st(n,K_{1,k}) \leq \floor{n^2\!/4}$.
Equality is attained by the bipartite Tur\'an graph.

Similarly, if $k+1\leq n \leq  2k-2$ and $m > \floor{(k-1)n/2}$,
 an edge $vw$ satisfies $d(v)+d(w) \geq \ceil{4m/n} \geq 2k-1 > n$, hence
 also here a triangle $vwz$ occurs, and assuming $d(v)\geq d(w)$ we have $d(v)\geq k$.
This yields a $K_{1,k}$ with a further edge in the neighborhood of $v$,
 implying $\st(n,K_{1,k}) \leq \floor{(k-1)n/2}$.
Equality is attained by any $(k-1)$-regular graph of order $n$ if
 $(k-1)n$ is even, or otherwise by any graph with one vertex of degree $k-2$
 and all the other $n-1$ vertices of degree $k-1$.
\epf

\btm
\label{t:xx}
Let $t \geq  2$ and $n \geq  2t$.
 \begin{itemize}
  \item[$(i)$]  If $t = 2$, then
  $$
      \st(n,2K_2) =
 \begin{cases}
  \begin{tabular}{lll}
    $n - 1 = \ex(n, P_4) = \ex(n, 2K_2)$ && $\mathrm{if } \ 3 \nmid n \,,$ \\
    $n = \ex(n, P_4) +1 = \ex(n, 2K_2) +1$ && $\mathrm{if } \ 3 \mid n \,.$
  \end{tabular}
 \end{cases}
  $$
  \item[$(ii)$]  If $t \geq  3$, then $\st(n,tK_2) = \ex(n,tK_2)$.
 \end{itemize}
\etm
 
\bpf
Since $D(2K_2) = \{P_4\}$, by Lemma \ref{l:equiv} we have $\st(n,2K_2) = \ex(n,P_4)$, and the rest  for $(i)$ follows from the known values of $\ex(n,P_4)$ and $\ex(n,2K_2)$.
 
Concerning $(ii)$ we clearly have $\ex(n,tK_2) \leq \st(n,tK_2)$.
Consider now a graph $G$ with $n \geq  2t$ vertices and $\ex(n,tK_2) +1$ edges. 
Assume for a contradiction that $G$ contains no proper supergraph of $tK_2$.

Since $e(G)>\ex(n,tK_2)$, $G$ contains $tK_2$ and it must be induced by assumption.
Let $M$ be the largest induced matching in $G$, then  $|M| \geq  t$.
Consider the edges not in $M$.
Every such edge must be incident with precisely one $M$-edge, hence one $M$-vertex:
 at most one as $M$ is induced, and at least one as $M$ is largest and none of its edges is extendible to a supergraph of $2K_2$ with an external edge.

If $v$ is a vertex in $V(G) -M$ and $v$ is adjacent to two vertices of $M$, then those vertices must belong to the same edge $e$ in $M$, else replacing the two neighbor edges with a $P_4$ containing $v$ a supergraph of $2K_2$ would occur.
It follows that, together with the single neighbor edge, $v \cup e$ induces a $K_3$ which must be isolated, as any edge incident with this $K_3$ would create $P_4$ and force a non-induced $tK_2$.

Further, there are no two independent edges incident to distinct vertices of the same edge in $M$, as they would form $P_4$ and then $M$ would contain a non-induced $tK_2$.
Hence the edges of $G$ form vertex-disjoint stars or copies of $K_3$,  but then  $e(G) \leq  n  < 2n-3 \leq \ex(n,tK_2)$ for $t \geq 3$, a contradiction. 
Consequently, for $t \geq  3$ we have $\st(n,tK_2) = \ex(n,tK_2)$.
\epf

\btm
\  $\st(n,P_4) = \ex(n,\{C_3 ,C_4\})$.
\etm

\bpf

First observe $\ex(n+3  ,\{C_3,C_4\}) \geq  \ex(n ,\{C_3 ,C_4 \}) + 3$ holds for all $n\geq 1$. 
Indeed, suppose $H$ realizes  $\ex(n,\{C_3 ,C_4 \})$.
Add three vertices  $x ,y ,z$ forming a path $P_3$ and take one new edge connecting $H$ to this $P_3$ to get $H^*$ which is still $C_3$-free and $C_4$-free. 

Now clearly $\ex(n,\{C_3 ,C_4\}) \leq  \st(n,P_4)$ because if there is no $C_3$ and no $C_4$ then even if there is a $P_4$ it must be induced. 

Consider now any graph $H$ with $n$ vertices and $\ex(n,\{C_3 ,C_4 \}) +1$ edges.  
If $H$ contains $C_4$, we are done as it is on four vertices and strictly contains $P_4$.
So, we may assume that $H$ is $C_4$-free, hence contains $C_3$.  
If this $C_3$ has an edge incident with it, we are done as this forms a graph on 4 vertices strictly containing $P_4$.
Consequently every triangle in $H$ is a connected component.

Removing an isolated  $C_3$ we obtain a graph $H_1$ with $n-3$ vertices  and $(\ex(n-3 , \{C_3 ,C_4 \}) +1)-3\geq \ex(n-3 , \{C_3 ,C_4 \}) +1$  edges. 
But it cannot contain $C_4$  as neither $H$ did, so it contains $C_3$ and the argument can be repeated.
In this way the structure is shrinked down by isolated vertex-disjoint copies of $C_3$   until we reach $n'$ = 1 or 2 or 3 vertices.

The cases $n'=1$ and $n'=2$ are impossible because those $(C_3,C_4)$-free extremal graphs are complete, so no graph with more edges can exist with those small orders.
The only possibility would be $n'=3$ where $\ex(3,\{C_3 , C_4 \}) = 2$ and one further edge yields $K_3$, thus originally $H=\frac{n}{3}C_3$ and
 before the last reduction step we would have $2C_3$, i.e., six vertices and six edges.
However $\ex(6, \{C_3 ,C_4\} ) \geq 6$ as demonstrated by $C_6$, so by assumption we should have at least $7>e(2C_3)$ edges before the last step.
This final contradiction completes the proof.  
\epf

We now proceed to estimate $\st(n,T)$ via the diameter of $T$.
For that, let us first state an auxiliary fact.

\blm
Let $T$ be a tree of order $q +1$.
If a graph $G$ has minimum degree at least $q$, then every subtree $T^*$ of $T$ can be completed to $T$ in $G$. 
\elm

\bpf
Suppose we get stuck at a subtree $T^{**}$ while trying to expand $T^*$ to $T$. Then there is a vertex $u$ in $T^{**}$ that has a neighbor $w$ in $T$ that was not added yet.  In particular, if $|T^{**}| \leq  q$,  then $u$ is adjacent to at most $q - 1$ vertices of $T^{**}$ but the degree of $u$ in $G$ is at least $q$, so we can add to it the neighbor $w$ from outside $T^{**}$, contradicting that the process got stuck at $T^{**}$.       
\epf

\btm
\label{t:treediam}
Let $T$ be a tree  with diameter $k \geq  2$  and order $q \geq  k+1   \geq  3$.  
Then $\ex(n, \{C_3,\dots,C_{k+1} \}) \leq   \st(n,T) \leq  \ex(n,\{C_3,\dots,C_{k+1}\}) +(q-1)n$.   
\etm

\bpf
Lower bound: Clearly if a graph contains no cycle from $\{C_3,\dots,C_{k+1}\}$, then if a copy of $T$ occurs, it must be induced. 

\msk

\nin
Upper bound: Let $G$ have $n$ vertices and $\ex(n, \{C_3,\dots,C_{k+1} \}  +  (q-1)n  +1$ edges. 
Repeatedly remove vertices of degree at most $q-1$ until no such vertex remains, and call the remainder graph $F$. 
This $F$ is not empty, since we removed at most $(q-1)n$ edges, fewer than $e(G)$.

Further, $F$ must contain a cycle from  $\{ C_3 ,\dots,C_{k+1} \}$, for otherwise we would have $e(F) \leq  \ex( |F| , \{ C_3,\dots,C_{k+1 } \}) \leq  \ex(n, \{C_3,\dots,C_{k+1} \})$  by monotonicity of Tur\'an numbers, and the contradiction $e(G) \leq  e(F) + (q-1)n  \leq  \ex(n,\{C_3,\dots,C_{k+1}\})   + (q-1)n  < e(G)$ would follow. 

Hence $F$ has minimum degree at least  $q$ and contains a cycle $C_ s$, where $3  \leq  s \leq  k +1$.  
Consider $T$.
As  $s \leq   k+1$ = diam$(T) +1$, there are two non-adjacent vertices $u ,v$ in $T$ at distance $s - 1$ apart, and connecting them by an edge creates a cycle $C_s$ with all vertices in $T$, forming a unicyclic graph $H\supsetneq T$. 

We will show that this cycle $C_s$ can be extended to a copy of $H$  in $F$, hence a non-induced copy of $T$. 
Let $u = v_1,\dots,v_s = v$   be the vertices of $C_s$, and delete the edge $(u,v)$ from $C_s$.  After this deletion still  all  vertices in $F$ have degree  at least $q-1$ and we can complete a copy of $T$  from the path $u = v_1,\dots,v_s  = v$  originally in $T$, 
by the above lemma. 
Having completed $T$ we add back the edge $(u,v)$  to get $H$, and the proof is done.
\epf

\bpn
If $G$ is an acyclic disconnected graph on $k$ vertices, then $\ex(n,G) \leq   \st(n,G) \leq  \ex(n,T) \leq  ( k-2)n/2$. 
In particular, $\st(n,G)$ is bounded by a linear function of $n$ if and only if $G$ is an acyclic disconnected graph.
\epn

\bpf
Clearly $\chi(G) = 2$ and $D(G)$ contains a forest or a tree $T$ on $k$ vertices, and the upper bound follows by the recently announced proof of the Erd\H os--S\'os Conjecture;
 see \cite{RS26} for a simplified version.

If $G$ contains a cycle $C_k$, then already $\st(n,G)  \geq   \ex(n,C_k)$ is a non-linear lower bound.
Also, if $G$ is a tree, then Theorem \ref{t:treediam} proved that $\st(n,G)$ is not linear.
Hence $G$ must be an acyclic disconnected graph.
\epf

\brm
In case $G$ contains isolated vertices, we refer to Section \ref{s:isol} for better bounds; and in case $G  = tK_2$ we refer to Theorem \ref{t:xx}.
\erm

\bpn
The value of $st(n,G)$ is constant if and only if $G$ is edgeless or $G = K_2 \cup tK_1$ for some $t \geq 2$.
\epn

\bpf
As we already noted, the edgeless case is trivial: $\st(n,E_k)=0$ holds for all $n\geq k$.

Concerning nonempty graphs, assume first $e(G) \geq 2$.
If $G = P_3$ then $\ex(n,P_3) = \floor{n/2}$ and if $G  =  2K_2$ then $\ex(n,2K_2) = n -1$ (for $n \geq 4$), and
 also for any larger $G$ at least a linear growth is valid for $\st(n,G)$
 by $\ex(n,G) \leq \st(n,G)$ and monotonicity of the Tur\'an function.
 
Hence it remains to consider $G = K_2 \cup tK_1$, $t \geq 1$.   
For $t = 1$ we have $D(K_2 \cup K_1)  = \{P_3\}$, therefore $\ex(n,K_2 \cup K_1) = \ex(n, P_3)  = \floor {n/2}$.
 
For $t \geq 2$ any graph with two edges contains either $2K_2$ or $P_3 \cup K_1$ and in both cases properly contains $K_2 \cup K_1$.
\epf

\section{Even cycles and theta graphs}
\label{s:theta}

Here we investigate $\st(n,G)$ for the case of cycles $G=C_\ell$ of even length at least 6.
Recall that $\st(n,C_4)  = \ex(n,K_3)$ is proved in Proposition \ref{p:st=ex} $(i)$.
In the present context, the so-called theta graphs will occur in a natural way.
A (generalized) theta graph $\theta_{k_1,\dots,k_m}$ has two vertices
 of degree $m$ connected by $m$ internally disjoint paths of
 respective lengths $k_1,\dots,k_m$.
We shall use the particular case $m=3$ of the following extremal result.

\igs

\theos{D.}{Liu, Yang \cite{LY23}}{Let $k_1,\dots,k_m$ be positive integers with the same parity, in which $1$ appears at most once.
Then $\ex(n, \theta_{k_1,\dots,k_m}) = O(n^{1+1/k^*})$,
where $k^* = \frac{1}{2} \min _{1\leq i  < j\leq m} (k_i + k_j)$.}

\igs

We prove the following asymptotic estimates for even cycles.
Up to length $\ell=24$ the estimates are exhibited in Table \ref{tab:cyc}.

\btm
\label{t:cyc}
Let $\ell\geq 6$ be an even integer.
 \begin{itemize}
  \item[$(i)$] We have $\st(n,C_\ell) = O(n^{1+1/\ceil{\ell/4}})$
   for every $\ell$.
  \item[$(ii)$] If $\ell\in\{6,8,10,12,18,20\}$ then
   $\st(n,C_\ell) = \Theta(n^{1+1/\ceil{\ell/4}})$.
  \item[$(iii)$] If $\ell=14,16$ then
   $c_1n^{6/5} \leq \st(n,C_\ell) \leq c_2n^{5/4}$
    for some positive constants $c_1,c_2$.
  \item[$(iv)$] If $\ell > 20$ with $\ell=4k-2$ or $\ell=4k$, then
    for some constant $c_k > 0$ we have
    $$
      \st(n,C_\ell) \geq
  \begin{cases}
    \begin{tabular}{ccl}
      $c_k n^{1+\frac{2}{6k-6}}$ && $\mathrm{if } \ \ell=4k-2$\,,\\ 
      $c_k n^{1+\frac{2}{6k-2}}$ && $\mathrm{if } \ \ell=4k$\,.
    \end{tabular}
  \end{cases}
    $$
 \end{itemize}
\etm

\begin{table}
 \begin{center}
  \begin{tabular}{ccccccc}
   $\ell$-cycle && LB && $\st(n,C_{\ell})$ && UB \\
      \hline
   $C_6$ && $n^{3/2}$ && $\Theta(n^{3/2})$ && $n^{3/2}$ \\
   $C_8$ && $n^{3/2}$ && $\Theta(n^{3/2})$ && $n^{3/2}$ \\
   $C_{10}$ && $n^{4/3}$ && $\Theta(n^{4/3})$ && $n^{4/3}$ \\
   $C_{12}$ && $n^{4/3}$ && $\Theta(n^{4/3})$ && $n^{4/3}$ \\
   $C_{14}$ && $n^{6/5}$ && ? && $n^{5/4}$ \\
   $C_{16}$ && $n^{6/5}$ && ? && $n^{5/4}$ \\
   $C_{18}$ && $n^{6/5}$ && $\Theta(n^{6/5})$ && $n^{6/5}$ \\
   $C_{20}$ && $n^{6/5}$ && $\Theta(n^{6/5})$ && $n^{6/5}$ \\
   $C_{22}$ && $n^{16/15}$ && ? && $n^{13/12}$ \\
   $C_{24}$ && $n^{18/17}$ && ? && $n^{13/12}$
  \end{tabular}
   \caption{Asymptotic estimates on $\st(n,C_{\ell})$ for even lengths $6\leq \ell\leq 24$;
    growth order of lower bounds (LB) and upper bounds (UB) proved in Theorem \ref{t:cyc}. 
    \label{tab:cyc}}
 \end{center}
\end{table}

\bpf
Assume $\ell=2k$.
Consider $H:=\theta_{k,k,1}$ if $k$ is odd, or $H:=\theta_{k+1,k-1,1}$
 if $k$ is even.
Due to the parity distinction, $H$ is bipartite and its
 shortest cycle has length $k+1$ or $k$, respectively.
Observe that $\frac{k+1}{2} = \frac{\ell+2}{4} = \ceil{\frac{\ell}{4}}$
 if $k$ is odd, and $\frac{k}{2} = \frac{\ell}{4} = \ceil{\frac{\ell}{4}}$
 if $k$ is even.

This $H$ properly contains $C_\ell$ and satisfies the parity condition
 of Theorem~D, consequently
  $$
    \st(n,C_\ell) \leq \ex(n,H) = O(n^{1+1/\ceil{\ell/4}})
  $$
 as claimed in $(i)$.
The upper bounds for $(ii)$ and $(iii)$ follow as particular cases.

For a lower bound we note that $C_\ell$ with any chord contains a cycle
 of length at most $\ell/2+1$.
So, every graph of girth at least $\ell/2+2$ is a feasible construction
 for $\st(n,C_\ell)$, that is
  $$
    \st(n,C_\ell) \geq \ex(n,\{C_3,\dots,C_{\ell/2+1}\}) \,.
  $$
It suffices to take $C_4$-free bipartite graphs for $\ell=6$ and $\ell=8$,
 yielding the lower bound $\Theta(n^{3/2})$ from the projective plane
 incidence graph \cite{KST}.
The estimates $\Theta(n^{4/3})$ and $\Theta(n^{6/5})$ for
 $\ell=10,12,18,20$ follow from the constructions of Lazebnik et al.\
 \cite{LUW1} on graphs without $C_6$ and $C_{10}$, respectively,
 which imply tight growth orders of the Tur\'an numbers of these cycles.

Since the growth rate of $\ex(n,C_8)$ is not known, concerning
 $\st(n,C_{14})$ and $\st(n,C_{16})$ the best we can do for $(iii)$ is take
 the construction of girth 10 or 12, as in case of $\ell=18,20$.

For $C_\ell$ with $\ell>20$ we apply the curerntly best known lower bounds \cite{LUW2}
 on $\ex(n,\{C_3,\dots,C_{2s+1}\})$, which are
 $\Theta(n^{1+\frac{2}{3s-3}})$ for $s$ odd, and
 $\Theta(n^{1+\frac{2}{3s-2}})$ for $s$ even.
These imply the validity of $(iv)$.
\epf

\def \tsk {\theta_{s*k}}
\def \thk {\theta_{s*3}}
\def \tok {\theta_{s*5}}
\def \tkk {\theta_{s*9}}
\def \tssk {\theta_{1,s*k}}
 
\brm
\label{r:tsk}
Let $s\geq 2$, and let $k\geq 3$ be odd.
Consider the graph $G=\tsk$ in which two nonadjacent vertices $u,v$ are
 connected by $s$ internally disjoint paths of length $k$.
Inserting the edge $uv$ we get the graph $G^* = \tssk$.
Then, with the notation of Theorem D we have $k^* = (k+1)/2$
and obtain
$\ex(n,G^*) =  O(n ^{1 + 2/(k+1)})$. So, by Lemma \ref{l:equiv}, $\st(n,G) \leq  \ex(n,G^*)$. 

A lower bound comes from the fact that any edge we add to $G$ forces a cycle of length at most $k+1$, hence a lower-bound is $\ex(n, \{C_3,\dots,C_{k+1}\})$ as before.
\erm

This principle leads to the following result.

\btm
For every fixed $s\geq 2$ we have $\st(n,\thk) = \Theta(n ^{3/2})$,
 $\st(n,\tok) = \Theta(n ^{4/3})$, and $\st(n,\tkk) = \Theta(n ^{6/5})$.
\etm

\bpf
As sketched in Remark \ref{r:tsk}, upper bounds on the growth orders
 can be derived from Theorem D.
The matching lower bounds hold due to the known estimates on
 $\ex(n,C_\ell)$ for $\ell=4,6,10$.
\epf

\section{Bipartite graphs and a theorem of Erd\H os}
\label{s:bip}

Erd\H os \cite{Er65} proved in 1965 the following result.

\igs

\theos{E.}{Erd\H os, \cite[Theorem 2]{Er65}}{For any integers $p\geq 2$ and $q\geq 1$ we have $\ex(n,K_{p+q,p+q}-K_{q,q})=O(n^{2-1/p})$.}

\igs  

An immediate consequence of this result with $q=1$ is:

\bpn
\label{p:bip}
~

\begin{itemize}
 \item[$(i)$] If $G$ is a bipartite graph containing $K_{p,p}$, with both sides of order  $p+1$ and   $p^2 \leq  e(G)  \leq  (p+1)^2 -2$,
 then $\ex(n,K_{p,p}) \leq \st(n,G) \leq \ex(n,K_{p+1,p+1}-e) \leq  O(n^{2-1/p})$. 
 \item[$(ii)$] For any $p\geq q\geq 2$, if
 $G  =  K_{p,q} - e$  then  $\ex(n ,K_{p,q})/2  \leq  \st(n , G) \leq  O(n^{2-1/p})$.
\end{itemize}
\epn

\bpf
Since $G$ is missing at least two edges of being $K_{p+1,p+1}$, it is obvious that $G$ is strictly  contained in $K_{p+1,p+1} - e$; and on the other hand $G$ contains $K_{p,p}$.  
Hence by monotonicity of Tur\'an numbers and by Lemma~\ref{l:equiv}, $\ex(n,K_{p,p}) \leq  \ex(n,G) \leq  \st(n,G) \leq  \ex(n, K_{p+1,p+1} - e)  \leq  O(n^{2-1/p})$ by Theorem E. 
This proves $(i)$ and also the upper bound in $(ii)$.
For the lower bound in the latter, it suffices to note the well-known fact
 that every graph contains a bipartite subgraph with at least half
 of its edges; and apply this to a graph extremal for $K_{p,q}$.
\epf

\brm
When an asymptotic estimate of order $\Theta(n^{2-1/p})$ is known for $\ex(n, K_{p,p})$, we get the correct order of exponent of $n$ for $\st(n,G)$.  This is known for $p = 2$ and $p=3$  with exponents $3/2$ and $5/3$, respectively.  
\erm

\bex
Consider $G=K_{3,3}-e$.
Then, by Lemma \ref{l:equiv} and on applying Brown's theorem \cite{Br66}, $\st(n,G)  = \ex(n , \{K_{3,3}$ \emph{and some 3-chromatic graphs containing} $K_3 \}) \leq \ex(n,K_{3,3}) = \frac{1}{2} n^{5/3}  +o(n^{5/3})$.
On the other hand,
suppose  $H(n)$  is an extremal graph for $K_{3,3}$.
We recall that there is a cut $( A ,B)$ of $V(H)$ containing at least half of the edges of $H$. Consider the bipartite graph $H^*$ with only the edges of the cut.  
This $H^*$ contains at least $\frac{1}{4} n^{5/3}  +o(n^{5/3})$ edges and no subgraphs of chromatic number at least  3,  and no copy of $K_{3,3}$,  as $H$ doesn't either.  
Hence $\st(n,K_{3,3} - e) = \Theta (n^{5/3})$. 
\eex

\brm
A further consequence of Theorem E above is as follows.

Denote by $H(p,l,l )$ the graph obtained from $K_{p,p}$ by adding in both sides $l-p$ vertices, each of them adjacent to all the $p$ vertices of $K_{p,p}$ in the other side. 
Erd\H os proved  $\ex(n, H(p, l, l) )  \leq  c(p,l)\,n^{2-1/p}$, for some constant $c(p,l)$.

Considering $\st(n, H(p,l,l ))$, the family $D(H(p,l,l ))$ contains  graphs of chromatic number $3$ by adding an edge in one side, and otherwise no matter how we add the edge we get a graph containing $K_{p+1,p+1}$.  
So,  using the cut-trick applied in the proof of Proposition \ref{p:bip} we dispose of the members of $D(H(p,l,l))$ of chromatic number~$3$,
 only with those containing $K_{p+1,p+1}$ remain (in fact $H(p+1,l-1,l-1 )$ only); hence $\st(n, H(p,l,l) )  \geq \ex(n,K_{p+1,p+1})/2$. 

Using again Theorem E for  $H(p+1,l-1,l-1)$  and the right side of
Lemma \ref{l:equiv}, we get:
$\st(n,H(p,l ,l ) ) \leq \ex(n, D(H(p +1 ,l-1, l-1 )) )   \leq  c(p+1, l-1)\,n^{2- 1/(p+1)}$.
Hence whenever $\ex(n,K_{p+1,p+1})  \sim  n^{2-1/(p+1)}$,  the sharp  exponent for $\st(n , H(p,l.l))$ is obtained. 
\erm

\bex
~

\nin
$(i)$\quad If $p = 1$, then $H(1,l ,l )$ is the double star  $S_{l,l}$ and we get $cn^{3/2} \leq  \ex(n, K_{2,2})/2 \leq  \st(n,S_{l, l}) \leq  c(2,l-1)\,n^{3/2}$.
This growth rate is valid for all double stars. 
(This comes also from the case of trees vs.\ diameter in
 Section~\ref{s:trees}, as the diameter of double stars is $3$.)

\nin
$(ii)$\quad If $p = 2$, then  $cn^{5/3}  \leq  \ex(n, K_{3,3} )/2  \leq   \st(n, H(2,l,l)) \leq  c(3,l-1)\,n^{5/3}$.
\eex

\begin{table}
\begin{center}
\renewcommand{\arraystretch}{1.3}
\begin{tabular}{ccccccc}
$G$  &  & $\ex(n,G)$ &  & $\st(n,G)$ &  & $D(G)$ \\
\hline
$2K_1$  &  & $-$ &  & 0 &  & $K_2$ \\
$K_2$  &  & 0 &  & $\binom{n}{2}$ &  & $-$ \\
$3K_1$  &  & $-$ &  & 0 &  & $K_2\cup K_1$ \\
$K_2\cup K_1$  &  & 0 &  & $\floor{\frac{n}{2}}$ &  & $P_3$ \\
$P_3$  &  & $\floor{\frac{n}{2}}$ &  & $\floor{\frac{n^2}{4}}$ &  & $K_3$ \\
$K_3$  &  & $\floor{\frac{n^2}{4}}$ &  & $\binom{n}{2}$ &  & $-$ \\
$4K_1$  &  & $-$ &  & 0 &  & $K_2\cup 2K_1$ \\
$K_2\cup 2K_1$  &  & 0 &  & 1 &  & $\{P_3\cup K_1,2K_2\}$ \\
$P_3\cup K_1$  &  & $\floor{\frac{n}{2}}$ &  & $\floor{\frac{2n}{3}}$ &  & $\{P_4,K_3\cup K_1,K_{1,3}\}$ \\
$K_3\cup K_1$  & & $\floor{\frac{n^2}{4}}$ &  & $\floor{\frac{n^2}{4}}$ &  & $K_4-P_3$ \\
$2K_2$  &  & $n-1$ &  & $n-\epsilon_n$ &  & $P_4$ \\
$P_4$  &  & $n-\epsilon_n$ &  & $\ex(n,\{C_3,C_4\})=\Theta (n^{3/2})$ &  & $\{K_4-P_3,C_4\}$ \\
$C_4$  &  & $\frac{1}{2}n^{3/2}+O(n)$ &  & $\floor{\frac{n^2}{4}}$ &  & $K_4-e$ \\
$K_{1,3}$  &  & $n$ &  & $\floor{\frac{n^2}{4}}$ &  & $K_4-P_3$ \\
$K_4-P_3$  &  & $\floor{\frac{n^2}{4}}$ &  & $\floor{\frac{n^2}{4}}$ &  & $K_4-e$ \\
$K_4-e$  &  & $\floor{\frac{n^2}{4}}$ &  & $\floor{\frac{n^2}{3}}$ &  & $K_4$ \\
$K_4$  &  & $\floor{\frac{n^2}{3}}$ &  & $\binom{n}{2}$ &  & $-$ 
\end{tabular}
\caption{Comparison of $\ex(n,G)$ and $\st(n,G)$ for graphs $G$ of order at most 4 and $n\geq |G|$.
Recall that $D(G)$ consists of the graphs $H=G+e$ ($|H|=|G|$) and $\st(n,G)  = \ex(n,D(G))$ where the last expression is insensitive to isolated vertices.
For $2K_2$ and $P_4$ we have $\epsilon_n=0$ if $3\,|\,n$ and $\epsilon_n=1$ otherwise.}
\end{center}
\end{table}

\section{Concluding remarks and open problems}
\label{s:open}

The strong Tur\'an parameter $\st(n,G)$ offers a structure-sensitive extremal problem and many open questions come to  mind.
 
The basic relation $\st(n,G) = \ex( n, D(G))$ connects the strong Tur\'an numbers to the classical Tur\'an numbers via the well-structured family $D(G)$ and, as shown, is calling for applications of classical extremal graph theory results.
 
An interesting general problem for further research is to find estimates on $\st(n, G\,\Box\,H)$  where `$\,\Box\,$' stands for the Cartesian product (also called box product) of two graphs.
 
A concrete case of this general problem is the following.

\bpm
\label{pm:grid}
Is it true that $st(n, P_p\,\Box\,P_q) = \Theta(n^{3/2})$ for every $p \geq 3$ and $q \geq 2$\,?
\epm

\brm
It is known \cite{BJST} that $\ex(n, P_t\,\Box\,P_t) = \Theta(n^{3/2})$ for every $t \geq 2$.
However, the asymptotic asked in Problem \ref{pm:grid}
 is not valid for $p=q=2$, as $P_2\,\Box\,P_2 \cong C_4$ and we have $\st(n,C_4) = \ex(n,K_3) = \Theta(n^2)$.
On the other hand it is valid for $p=3$, $q=2$ as $D(P_3\,\Box\,P_2)$
 consists of $K_{3,3}-e$ and some $3$-chromatic graphs, and
  $\ex(n,K_{3,3}-e) = \Theta(n^{3/2})$ by a theorem of Erd\H os (Theorem E above).
\erm

\paragraph{Acknowledgement.}
Research of the second author was supported in part
by the ERC Advanced Grant no.\ 101054936 ``ERMiD''.

\end{document}